\documentclass{article}
\usepackage{amsmath,amssymb}

\newtheorem{prop}{Proposition}
\newtheorem{theo}[prop]{Theorem}
\newtheorem{lemm}[prop]{Lemma}

\newcommand{\defi}{\addtocounter{prop}{1}\par\noindent\textbf{Definition \arabic{prop}.} }
\def\proof{\noindent\textit{Proof. }}
\newcounter{example}

\def\C{{\mathbb C}}

\def\I{\mathbf{I}}

\def\qed{\ \hfill\raise3.5pt\hbox{\framebox[1ex]{\ }}}

\newcommand{\X}[2]{X_#1^{(#2)}}
\newcommand{\p}[2]{P_#1^{(#2)}}

\begin{document}
%%%%%%%%%%%%%%%%%%%%%%%%%%%%%%%%%%%%%%%%%%

\title{Spectral  parameter power series
for arbitrary order linear differential equations}
\author{Vladislav V.\ Kravchenko, R.\ Michael Porter, Sergii M.\ Torba \\
Department of Mathematics, Cinvestav-Quer\'etaro, Mexico}

\maketitle

\begin{abstract}
 Let $L$ be the $n$-th order linear differential operator
$Ly = \phi_0y^{(n)} + \phi_1y^{(n-1)} + \cdots + \phi_ny$ with variable coefficients.
 A representation is given for $n$ linearly independent solutions
of $Ly=\lambda r y$ as  power series in $\lambda$, generalizing
the SPPS (spectral parameter power series) solution which has been
previously developed for $n=2$. The coefficient functions in these
series are obtained by recursively iterating a simple integration
process, begining with a solution system for $\lambda=0$.  It is shown
how to obtain such an initializing system working upwards from equations
of lower order. The values of the successive derivatives of the power
series solutions at the basepoint of integration are given, which
provides a technique for numerical solution of $n$-th order initial
value problems and spectral problems.
\end{abstract}



\section{Introduction}

The spectral parameter power series (SPPS) representation for
solutions of second-order linear differential equations, introduced in
\cite{KrCV08} and developed into a numerical technique in
\cite{KrPorter2010}, has proved to be an efficient tool for solving
and studying a variety of problems (see the review \cite{KKRosu}). The
SPPS method offers a simple and numerically stable procedure for
computing coefficients of the Taylor series of the solution with
respect to the spectral parameter. The procedure consists of recurrent
integration and requires knowledge of a nonvanishing solution of the
equation with one fixed value of the spectral parameter (which can be
zero). The availability of such a solution leads to a Polya
factorization of the second-order differential operator and to a
convenient form for its right inverse operator.  In \cite{KKB} the
SPPS representation was obtained for fourth-order linear differential
equations of a special form.

Here we extend these ideas to linear differential equations of
arbitrary order. In contrast to the second order case, the
existence of a single nonvanishing solution is insufficient for obtaining a
convenient Polya factorization of the differential operator. Instead,
for factoring a linear differential operator of order $n$ a system
of $n$ linearly independent solutions (for a fixed value of the
spectral parameter) is required whose $n$ corresponding partial
Wronskians do not vanish. In \cite{CdS} the existence of such a system
of solutions was proved, and in the case of real valued coefficients
there is even an algorithm for its construction, see \cite{Mammana}.

We describe a procedure for computing the coefficients of the SPPS
representation of a generating set of solutions of the $n$-th order
linear differential equation, given a system of solutions leading to a
Polya factorization and using a simple procedure of recurrent
integration. We note that different procedures for constructing the
SPPS coefficients were considered also in \cite{Fage} and
\cite{Leontiev}. The Polya factorization and corresponding systems of
solutions and iterated integrals are introduced in Section
\ref{sec:systems}. In Section \ref{sec:solution} the SPPS
representation for solutions of linear $n$-th order differential
equations is defined and justified, and an example is worked out in
Section \ref{sec:example}.

\section{Systems related to a general solution of the homogeneous
  equation\label{sec:systems}}

 \subsection{P\'olya factorization\label{subsec:polfac}}

Let $n\ge2$, and let $L$ be the general linear ordinary differential operator
\begin{equation}  \label{eq:defL}
  Ly = \phi_0y^{(n)} + \phi_1y^{(n-1)} + \cdots + \phi_ny
\end{equation}
on a real interval $x_1<x<x_2$.  We will always assume that the
$\phi_j$ are continuous and that $\phi_0\equiv1$. Write
$\partial= \partial/\partial x$. The P\'olya factorization \cite{Pol}
\begin{equation}  \label{eq:polfac0}
  L = \frac{W_n}{W_{n-1}}\partial  \frac{W_{n-1}^2}{W_{n-2}W_n}\partial
\cdots \partial \frac{W_2^2}{W_1W_3}\partial \frac{W_1^2}{W_0W_2}\partial
 \frac{1}{W_1}
\end{equation}
of $L$ can be constructed by taking linearly independent solutions
$y_1,y_2,\dots,y_n$ of $Ly=0$, satisfying the condition that the
following $n+1$ functions do not vanish:
\begin{equation}\label{eq:W}
W_0=1,\ W_1=y_1,\ W_2=W[y_1,y_2],\ \dots,\
   W_n = W[y_1,y_2,\dots,y_n].
\end{equation}
  Here
 \[ W_j = \left| \begin{array}{ccccc}
  y_1&y_2&\dots&y_j\\
  y_1'&y_2'&\dots&y_j'\\
  \vdots& \vdots &\vdots &\vdots\\
  y_1^{(j-1)}&y_2^{(j-1)}&\dots&y_j^{(j-1)}
          \end{array} \right|
\]
is the Wronskian operator.  The existence of such a solution system
$(y_k)$ is established in \cite{CdS}, where it can be seen that in fact almost all complex-valued solution sets satisfy this nonvanishing requirement.  We will write
\eqref{eq:polfac0} for convenience as
\begin{equation}  \label{eq:polfac1}
 Ly = \frac{1}{b_n}\partial \frac{1}{b_{n-1}}\partial  \cdots
   \frac{1}{b_2}\partial
  \frac{1}{b_1}\partial \frac{1}{b_0} y,
\end{equation}
noting that $b_0,\dots,b_n$ are nonvanishing
and
\begin{equation}  \label{eq:Lb0}
   Lb_0 =0.
\end{equation}

\subsection{Derivatives of iterated integrals and the P\'olya system\label{subsec:derivii}}

We will further simplify the notation by writing $\int\!\!f$ to
signify the indefinite integral $\int_{x_0}^{(\cdot)}f(s)\,ds$ from a
fixed basepoint $x_0$, where $x_1<x_0<x_2$. In addition, we introduce
the collection $\I_{i,j}$ of linear operators defined on integrable
functions $y$ by
\[ \I_{i,j}[y] = \int b_i\int b_{i+1}\cdots \int b_{j-1}\int b_j y
\]
when $i\le j$. In particular $\I_{j,j}y=\int b_jy$, and additionally
we specify
\[ \I_{i,i-1}[y] \equiv 1, \]
and that $\I_{i,j}[y] =0$ for $j\le i-2$.  Following the common usage
of a symbol to denote a function as well as the operation of
multiplication by that function, we will freely write
\[    \I_{i,j} =  \I_{i,j}[1] =
     \int b_i\int b_{i+1}\cdots \int b_{j-1}\int b_j
   \]
without fear of confusion.
Clearly
\begin{equation}   \label{eq:recurrence}
 \I_{i,j}' = b_i \I_{i+1,j}
\end{equation}
is valid for all values of $i,j\ge0$.

The \textit{P\'olya system} \cite{dVL} determined by the functions $b_0,\dots,b_n$  is the collection of functions
\begin{equation}     \label{eq:polSys}
b_0\I_{1,k} \quad\quad (0\le k\le n-1) .
\end{equation}

The following proposition states that the P\'olya system \eqref{eq:polSys} represents a fundamental system of solutions for the equation $Ly=0$ where the differential operator  $L$ is being given by \eqref{eq:polfac1}.

\begin{prop}[{\cite[Chap. 3, \S 11]{dVL}}]\label{prop:basicbasis}
  The functions $b_0\I_{1,k}$ $(0\le k\le n-1)$ forming the P\'olya
  system are linearly independent, and satisfy $L(b_0\I_{1,k}) = 0$.
\end{prop}

In relation with the P\'olya system let us consider the $\ell$-th derivatives of the functions $b_0\I_{i,j}$. From the relation \eqref{eq:Ideriv} the following is easily verified by induction.

\begin{lemm}\label{lemm:basis} Let $\ell\ge0$. The $\ell$-th derivative
of $b_0 \I_{i,j} $ can be expressed as
  \begin{equation}\label{eq:Ideriv}
    \partial^\ell \bigl(b_0\I_{i,j}\bigr) = \sum_{\alpha=0}^{\ell}
    A_{\ell,\alpha}\I_{i+\alpha,j},
  \end{equation}
where the coefficients $A_{\ell,\alpha}$ $(0\le\alpha\le\ell)$ are
  polynomials of degree $\alpha+1$ in the functions
  $b_0,b_i,b_{i+1},\dots,b_{i+\alpha-1}$ and their derivatives, with integer
  coefficients which do not depend on $i,j$.  These polynomials are
  defined recursively starting with $A_{0,0}= b_0$, and then for $\ell \ge 1$
\begin{equation}\label{eq:Aellalpha}
A_{\ell,\alpha}=
    \begin{cases}
      A_{\ell-1,\alpha}', & \alpha=0,\\
      A_{\ell-1,\alpha}' + A_{\ell-1,\alpha-1}b_{i+\alpha-1},& 1\le\alpha\le\ell-1,\\
       A_{\ell-1,\ell-1}b_{i+\ell-1},& \alpha=\ell.
    \end{cases}
\end{equation}
More precisely, each expression $A_{\ell,\alpha}$ is a linear
combination of expressions
$\partial^{\ell_0} b_0 \partial^{\ell_1} b_i
 \partial^{\ell_2} b_{i+1}\ldots
 \partial^{\ell_{\alpha}} b_{i+\alpha-1}$
with $\ell_k\ge 0$ and
$\ell_0+\ell_1+\ldots+\ell_{\alpha}=\ell-\alpha$.  The final
coefficient $A_{\ell,\ell}$ is nonvanishing.
\end{lemm}

Note that \eqref{eq:Aellalpha} may contain functions $b_{i+k}$ for
which $i+k>n$. These $b_{i+k}$ can be taken to be arbitrary
nonvanishing since the final expression is operated on by
$\I_{i+\alpha,j}\equiv 0$.

It is easily seen that \eqref{eq:Ideriv} holds even when $\I_{i,j}$ is replaced by
$\I_{i,j}f$ for any function $f$, so long as the order $\ell$ of the derivative
does not exceed $j-i+1$.

\section{Solution in terms of formal powers\label{sec:solution}}

\subsection {SPPS solutions of $n$-th order\label{subsec:SPPS}}

Apart from $\phi_1,\dots,\phi_n$, which have been assumed continuous on
the closed interval $[x_1,x_2]$, we consider an additional continuous
function $r$ defined on the same interval.  Define
\begin{equation} \label{eq:Pkm}
  \p{k}{m} = (mn+k-1)! \left(\I_{1,n} rb_0\right)^m  \I_{1,k-1}
\end{equation}
for $1\le k\le n$ and all $m\ge0$.  Repeated application of
\eqref{eq:polfac1}, which says that $L^{-1}=b_0\I_{1,n}$ is a right
inverse for $L$, gives us the differential relation
\begin{equation} \label{eq:LPkm}
  L[ b_0 \p{k}{m}] = \frac{(mn+k-1)!}{((m-1)n+k-1)!} rb_0 \p{k}{m-1}.
\end{equation}
We will refer to $\p{k}{m}$ as the \textit{main formal powers} for $L$ because
of the following.

\defi The \textit{SPPS solution system} $u_k(x)=u_k(x;\,\lambda)$ for
$L$ and $r$ is the collection of functions
\begin{equation}  \label{eq:uk}
   u_k(x)  = b_0\sum_{m=0}^\infty
 \frac{1}{(mn+k-1)!}  \p{k}{m}(x)\lambda^m
\end{equation}
defined on $[x_1,x_2]$ for $1\le k\le n$.  For $n=2$, the powers $\p{1}{m}$ and
$\p{2}{m}$ have been denoted   $X^{(2m+1)}$and   $\widetilde X^{(2m)}$ in most of the
SPPS literature for equations of second order, as in \cite{KKRosu,KMoT,KrPorter2010,Porter}.
Here for $1\le k\le n$ we define the $k$-th sequence of \textit{secondary
  formal powers} for $L$ and $r$ by setting  $\X{k}{j}\equiv 0$ for all $j<0$, $\X{k}{0}\equiv1$ and
then recursively for $j\ge 1$
\begin{equation}  \label{eq:defXkj}
\X{k}{j} =
\begin{cases}
 j\, \I_{k-j,k-j}\X{k}{j-1},& 1\le j\le k-1,\\
 j\, \I_{n,n}b_0r\X{k}{j-1},&   j\equiv k \mod n,\\
 j\, \I_{n-j',n-j'}\X{k}{j-1},& j\equiv k+j'  \mod n,\ 1\le j'\le n-1.
\end{cases}
\end{equation}
In loose terms, we multiply by $b_{k-j \mod n}$ and then integrate, except when
$j\equiv k$ mod $n$, in which case we multiply by $b_nb_0r$ before
integrating.  In the third case of \eqref{eq:defXkj}, clearly
$\partial \X{k}{j} =j\, b_{n-j'}\X{k}{j-1}$.  Further, by comparing with
\eqref{eq:Pkm} we see that the subsequence of main formal powers is
included as $\p{k}{m}=\X{k}{mn+k-1}$.

\begin{theo}\label{theo:SPPS}
  The SPPS series $u_k(x)$ ($1\le k\le n$) converge uniformly on the closed interval
  $[x_1,x_2]$ for each fixed $\lambda\in\C$.  They are linearly
  independent, and satisfy
\begin{equation}  \label{eq:Luk}
    Lu_k = \lambda r u_k.
  \end{equation}
  Further, for $1\le\ell\le n-1$, the $\ell$-th derivative of $u_k$ is
  given by substituting
\begin{equation}  \label{eq:dPkm}
  \partial^{\ell}\bigl(b_0\p{k}{m}\bigr) =
  \sum_{\alpha=0}^\ell (mn+k-\alpha)_\alpha A_{\ell,\alpha}\,\X{k}{mn+k-(\alpha+1)}
\end{equation}
in place of $b_0\p{k}{m}$ in \eqref{eq:uk}, where $(x)_n$ is the Pochhammer symbol and the coefficients
$A_{\ell,\alpha}$ are the functions defined in Lemma \ref{lemm:basis}
for the particular case $i=1$, $j=
\begin{cases}
k, & \text{if } m =0,\\
n, & \text{if } m>0.
\end{cases}$
\end{theo}

Note that $A_{\ell,\alpha}$ in \eqref{eq:dPkm} refer to
$A_{\ell,\alpha}(b_0, b_0',\ldots, b_1,b_1',\dots,b_{\alpha-1}^{(\ell-\alpha)})$ since by \eqref{eq:Pkm}, the
$\p{k}{m}$ are obtained by integrations of the form $\I_{1,j}\cdots$ where $i=1$.

\proof First we show the convergence of \eqref{eq:uk}. By continuity,
there is a uniform bound
\[ |b_0(x)b_n(x)r(x)|\le M_1 ,\quad  |b_j(x)|\le M_1  \quad  (1\le j\le n-1),
\]
on $[x_1,x_2]$. To simplify the estimates we write
$e_j(x)=(M_1|x-x_0|)^j$.  Thus whenever $|f(x)|\le  e_j(x)$, both
$\I_{k,k}f(x)$, $1\le k\le n-1$, and $\I_{n,n}b_0rf(x)$ are bounded by $e_{j+1}(x)/(j+1)$.
Thus by induction $|\p{k}{m}(x)|\le e_{mn+k-1}(x)$
since \eqref{eq:Pkm} involves $mn+k-1$ integrations.
From this estimate,
\begin{align*}
  \sum_{m=0}^\infty\left| \frac{ \p{k}{m}(x) }{(mn+k-1)!}\lambda^m \right|
  &\le \sum_{m=0}^\infty  \frac{e_{mn+k-1}(x)}{(mn+k-1)!} |\lambda|^m \\
  &< ( M_1|x-x_0|)^{k-1}\sum_{m=0}^\infty \frac{(M_1|x-x_0||\lambda|^{1/n})^{mn}}{(mn)!} \\
  &\le  ( M_1|x-x_0|)^{k-1}  \sum_{j=0}^\infty \frac{(M_1|x-x_0||\lambda|^{1/n})^{j}}{j!} \\
  &<   \infty
\end{align*}
since the last sum is simply $\exp M_1|(x-x_0)\lambda|^{1/n}$.
This shows that the series \eqref{eq:uk} for $u_k$ converges uniformly
on $[x_1,x_2]$. (Bounds for $\p{k}{m}$ were also given in \cite{Leontiev}.)

The finite sum \eqref{eq:dPkm} is obtained from the remark following
Lemma \ref{lemm:basis}. To see that one may differentiate
\eqref{eq:uk} term by term, we refine our estimates.

By induction it is seen that $|\X{k}{j}|\le e_j(x)$.  We can take a
finite bound $M_2$ such that
\[ A_{\ell,\alpha}(\dots \partial^{\ell_i} b_{i+\alpha} \dots) \le M_2 \]
for all combinations $0\le \alpha\le \ell$ and $0\le\ell_i\le \ell-\alpha$.  Then for
$1\le\ell\le n-1$ we note that $\X{k}{i}$ involves
precisely $i$ integrations, so
\begin{align*}
  |\partial^{\ell}\bigl(b_0\p{k}{m}\bigr)|
  &\le  \sum_{\alpha=0}^\ell (mn+k-\alpha)_\alpha |A_{\ell,\alpha}|\,|\X{k}{mn+k-(\alpha+1)}| \\
  &\le  \sum_{\alpha=0}^\ell  (mn+k-\alpha)_\alpha M_2(M_1|x_2-x_1|)^{mn+k-(\alpha+1)}.
%   \\
%  &\le \sum_{\alpha=0}^\ell (mn+k-\alpha)_\alpha  M_2^{mn+k-\alpha}.
%  \\
%  &<  \frac{M_2^{k+mn+1}}{M_2-1}.
  \end{align*}
(The estimate holds even though the specific functions
$A_{\ell,\alpha}$ vary according to the indices $k,m$.) Therefore
\[  \sum_{m=0}^\infty \left| \partial^{\ell}\frac{b_0\p{k}{m}}{(mn+k-1)!}\lambda^m \right|
  \le M_2\sum_{j=0}^\infty \frac{\bigl(M_1(x_2-x_1)\bigr)^{j}}{j! } \left|\lambda^{\left[\frac{j-k}{n}\right]+1}\right|
  \ < \ \infty.
\]
This justifies the assertion that the $n-1$ termwise differentiated series
obtained from \eqref{eq:uk} all converge. Now we are also justified in applying
operator $L$ termwise, and applying \eqref{eq:Lb0} and \eqref{eq:LPkm} we obtain
\begin{align*}
   L[u_k] =&  \sum_{j=0}^\infty \frac{1}{(mn+k-1)!} L[b_0 \p{k}{m}] \lambda^m \\
    =& \sum_{j=1}^\infty  \frac{1}{((m-1)n+k-1)!} r b_0\p{k}{m-1} \lambda^m \\
    =& \sum_{j=0}^\infty r b_0\p{k}{m} \lambda^{m+1} \\
    =& \lambda r u_k .
\end{align*}
The linear independence of the $u_k$ will be proved in Section
\ref{subsec:basepoint}, which will complete the proof. \qed

\subsection{Initial conditions\label{subsec:basepoint}}

In this section we look at the values of the SPPS solutions and their
derivatives at the integration base point $x_0$.  Most of the summands
in the formulas \eqref{eq:uk}, \eqref{eq:dPkm} vanish at
$x=x_0$. Indeed, since indefinite integrals always vanish at their
base point, $\X{k}{j}(x_0)=0$ vanishes unless $j=0$. From the fact
that $\p{k}{0}=(k-1)!\,\I_{1,k-1}$ it follows from \eqref{eq:uk} that
\begin{equation}  \label{eq:Pkx0 }
   u_1(x_0) = \begin{cases}
      b_0(x_0),  & k=1,\\
        0,        & k>1.  \end{cases}
\end{equation}

As to the derivatives of the $u_k$, the only terms which survive upon
evaluating \eqref{eq:dPkm} at $x=x_0$ are those for which the
exponents $k-(\alpha+1)+mn$ are equal to zero. Together with the
conditions $0\le\alpha\le\ell$, $1\le k\le n$, and $1\le\ell\le n-1$
this implies that $\partial^\ell\bigl(b_0\p{k}{m}\bigr)$ vanishes unless
$m=0$ and $\alpha=k-1$. Thus
\[ \partial^\ell\bigl(b_0\p{k}{m}\bigr)(x_0)=0 \qquad \text{ for } m>0, \]
and referring again to
\eqref{eq:dPkm}, we deduce that
\begin{equation}  \label{eq:dpk0x0}
  \partial^\ell\bigl(b_0\p{k}{0}\bigr)(x_0) = \begin{cases}
     0,  &  1\le \ell < k-1,      \\
     (k-1)!A_{\ell,k-1}(x_0), \ &  k-1\le\ell\le n-1.
          \end{cases}
\end{equation}
As a consequence, it follows that
\begin{equation}  \label{eq:dukx0}
     \partial^\ell u_k(x_0) =  \begin{cases}
   0 ,  &  \ell<k-1, \\
   A_{\ell,k-1}(x_0),
\ &  \ell \ge k-1 . \end{cases}
\end{equation}
for $1\le\ell\le n-1$ and $1\le k\le n$.
An important consequence of this formula is that
$\partial^\ell u_k(x_0)$ does not depend on $\lambda$: the coefficients
of $\lambda^m$ all vanish for $m\ge1$.  A further consequence of
\eqref{eq:dukx0} is that the Wronskian matrix of $u_1,\dots,u_n$ is
lower triangular and its determinant is equal to the product
\[ \prod_{k=1}^n \partial^{k-1} u_k(x_0) = \prod_{k=1}^n A_{k-1,k-1}(x_0)
\]
which is nonzero. This implies that these $n$ solutions of
$Lu-\lambda ru=0$ are linearly independent, and thus finishes the
proof of Theorem \ref{theo:SPPS}.

\subsection{On construction of a particular solution system}
The solution system $y_1,\ldots,y_n$ of $Ly=0$ can be constructed by applying Theorem \ref{theo:SPPS} inductively.

Indeed, consider an equation
\begin{equation}\label{eq:z}
    z^{(n-1)}+\phi_1 z^{(n-2)}+\ldots+\phi_{n-1} z = 0
\end{equation}
and let $z_1,\ldots,z_{n-1}$ be $n-1$ linearly independent solutions of \eqref{eq:z}. Then the functions
\[
\tilde z_0 = 1,\quad \tilde z_1= \int z_1,\quad \tilde z_{n-1}=\int z_{n-1}
\]
form a linearly independent family of solutions for the equation
\begin{equation}\label{eq:y0}
    y^{(n)}+\phi_1 y^{(n-1)}+\ldots+\phi_{n-1} y' = 0.
\end{equation}
Since almost all linear combinations with complex coefficients of $\tilde z_0,\ldots, \tilde z_{n-1}$ have non-vanishing Wronskians \eqref{eq:W} (we refer to \cite{CdS} for the precise meaning of the term ``almost all''), one may apply Theorem \ref{theo:SPPS} taking $r=\phi_n$ and $\lambda=-1$ to construct a complete family of solutions for the equation
\[
y^{(n)}+\phi_1 y^{(n-1)}+\ldots+\phi_{n-1} y' = -\phi_n y.
\]
And once again, taking a linear combination with complex coefficients one may easily find a solution system possessing non-vanishing Wronskians \eqref{eq:W}.

\section{Example\label{sec:example}}

We illustrate with a simple example the procedure presented in this
paper.  In more general instances it would be necessary to find
numerical approximations for the formal powers.

  Let $L=\partial^{n}$ and take $r(x)=1$. By inspection, the
functions
\[ y(x)\ =\ e^{\lambda^{1/n}x} \]
are linearly independent solutions of $Ly=\lambda y$ where
$\lambda^{1/n}$ varies over the complex $n$-th roots of 1. Let
$\lambda>0$. The above solutions can be written in complex form as
\begin{align*}
      e^{\lambda^{1/n}x} =& e^{\sqrt[n]{\lambda}\bigl(\cos\frac{2\pi k}{n}\bigr)x}   \left\{
     \cos\Bigl[\sqrt[n]{\lambda}\Bigl(\sin\frac{2\pi k}{n}\Bigr)x\Bigr] %\right. \\
 % &&\quad \left.
    +\ i \sin\Bigl[\sqrt[n]{\lambda}\Bigl(\sin\frac{2\pi k}{n}\Bigr)x\Bigr]\right\}.
\end{align*}
However, to obtain real-valued solutions from this approach involves
choosing the appropriate basis from among the $2n$ real and imaginary
parts, and is rather delicate, as it depends on the prime factorization
of $n$.

To apply Theorem \ref{theo:SPPS}, the basic solution set described in Subsection
\ref{subsec:polfac} is required, i.e.\  we begin with $\lambda=0$.
The homogeneous equation $y^{(n)}=0$ has the real solutions
\[  y_k(x) = \frac{x^{k-1}}{(k-1)!},\quad\quad  k=1,2,\dots,n. \]
The corresponding Wronskians are
\[  W_k = \begin{vmatrix}
  1 &x& \frac{x^{2}}{2!} &  \frac{x^{3}}{3!} & \cdots & \frac{x^{k-1}}{(k-1)!} \\
  0 & 1 &x& \frac{x^{2}}{2!}  & \cdots & \frac{x^{k-2}}{(k-2)!} \\
  0 & 0 & 1 &x& \cdots & \frac{x^{k-3}}{(k-3)!} \\
  &&&&\vdots&\\
  0 & 0 & 0 &0& \cdots & 1
                 \end{vmatrix}  = 1
\]
so $b_j=1$ for all $j$, and the P\'olya factorization amounts to the
obvious fact that $L=1\partial1\partial\cdots\partial1\partial1$.

Consequently, the construction of the SPPS formal powers and the SPPS
solutions for the equation $Ly=\lambda y$ is as follows. Taking as
base point $x_0=0$, we have
  \[ \I_{k-1} \ =\ \int1\cdots\int1 \ =\ \frac{x^{k-1}}{(k-1)!}, \]
and since $rb_0=1$, the main formal powers are
  \[ \p{k}{m}  = (mn+k-1)!(\I_{1,n}1)^m \I_{k-1}  =   x^{mn+k-1}, \]
a collection of ``true'' powers of $x$, indexed for the purpose at hand.
Formula \eqref{eq:uk}  yields the real valued solutions
\[ u_k(x) = \sum_{m=0}^\infty  \frac{x^{mn+k-1}}{(mn+k-1)!} \lambda^m
  \quad\quad(1\le k\le n)   \]
for the operator $\partial^n-\lambda$.

\section*{Acknowledgements}
Research of the first and third named authors was supported by CONACYT, Mexico via the project 222478.

\end{document}